\documentclass[12pt, a4paper]{article}

\usepackage{graphicx}

\usepackage[T1]{fontenc}
\usepackage{multirow}
\usepackage{float}
\usepackage[centertags]{amsmath}
\usepackage{amsfonts}
\usepackage{amssymb}
\usepackage{amsthm}
\usepackage{authblk}

\DeclareMathOperator{\Gal}{Gal}

\newtheorem{thm}{Theorem}[section]

\theoremstyle{definition}
 \newtheorem{defn}[thm]{Definition}
\theoremstyle{remark}
 \newtheorem{rem}[thm]{Remark}
\theoremstyle{example}
 \newtheorem{ex}[thm]{Example}
\theoremstyle{conjecture}
 \newtheorem{conj}[thm]{Conjecture}

\numberwithin{equation}{section}

\begin{document}

\begin{sloppy}

\title{Quantum groups: from Kulish-Reshetikhin discovery to classification}
\author[*]{Boris Kadets}
\author[*]{Eugene Karolinsky}
\author[**]{Iulia Pop}
\author[***]{Alexander Stolin}

\affil[*]{\small{Department of Mechanics and Mathematics, Kharkov National University}}
\affil[**]{Department of Mathematical Sciences, Chalmers University of Technology}
\affil[***]{ Department of Mathematics, Gothenburg University, Sweden}
\maketitle
\begin{abstract}
The aim of this paper is to provide an overview of the results about classification of quantum groups that were obtained in \cite{SP}, \cite{SP2}.


\end{abstract}

\section{Introduction}
The first example of a quantum group was found by Kulish and Reshetikhin in \cite{KR}. They discovered what was later named $U_{q}(\mathfrak{sl}_2)$ in relation to the study of the inverse quantum scattering method. Later, Drinfeld \cite{Dr} and Jimbo \cite{Ji} independently developed a general notion of quantum group. Today there are many different approaches to what quantum group is and the term has no clear meaning. Informally speaking, the quantum group is a deformation of a universal enveloping algebra of some Lie algebra $\mathfrak{g}$. Of course, the precise meaning should be given to the term deformation. We will be using the following definition.
\begin{defn}\label{Qgr}
 A quantum group is a topologically free cocommutative $\mathrm{mod}\ \hslash$  Hopf algebra over $\mathbb{C}[[\hslash]]$ such that $H/\hslash H$ is a universal enveloping algebra of some Lie algebra $\mathfrak{g}$ over $\mathbb{C}$.
\end{defn}
It is well known that many problems about Lie groups become simpler when they are written on the language of Lie algebras.
In general the existence of almost one-to-one correspondence between Lie groups and Lie algebras is one of the central parts of Lie theory.
Therefore,
it is desirable to obtain a notion of quantum algebra that will help to simplify problems about quantum groups gradually .
The first natural attempt was to look at the linear part of comultiplication of a quantum group $H$.
Indeed, one can define a co-Poisson structure $\delta: U(\mathfrak{g}) \to U(\mathfrak{g}) \otimes U(\mathfrak{g})$
by the formula $$\delta(x) = \frac{\Delta(a) - \Delta^{21}(a)}{\hslash}\  \mathrm{mod}\ \hslash, $$ where $x \equiv a\ \mathrm{mod}\ \hslash$.
Furthermore, from a co-Poisson structure on $U(\mathfrak{g})$ one gets a Lie bialgebra structure on $\mathfrak{g}$ and the
co-Poisson structure is uniquely determined by this Lie bialgebra structure.

The following problem naturally arises.

\begin{conj}[Drinfeld's quantization conjecture]
 Any Lie bialgebra can be quantized.
\end{conj}
Drinfeld conjecture was solved by Etingof and Kazhdan in \cite{EK1}, \cite{EK2}. However, it is not difficult to see that
one can get many different quantum groups from a Lie bialgebra (or, in informal terminology,
there are many different ways to quantize $(\mathfrak{g}, \delta)$).

Kazhdan and Etingof not only proved
Drinfeld quantization conjecture but found a correct
notion of quantum algebra. However, there can be many different quantum groups corresponding to a given Lie bialgebra (quantizations of this bialgebra).
They constructed a canonical co-Poisson structure on $U(\mathfrak{g}) \otimes \mathbb{C}[[\hslash]]$.
This structure is much finer then the co-Poisson structure discussed above.
The Lie groups -- Lie algebras correspondence has an analogy in the quantum world.

\begin{thm}\label{KazEti}
Let $\mathrm{Qgroup}$ be the category of quantum groups in the sence of Definition \ref{Qgr}.
Let $\mathrm{LieBialg}$ be the category of topologically free Lie bialgebras over
$\mathbb{C}[[\hslash]]$ with $\delta = 0\ \mathrm{mod}\ \hslash$. Then there exists a dequantization functor
$deQuant: \mathrm{Qgroup} \to \mathrm{LieBialg}$ that is an equivalence of categories.
\end{thm}

In their solution of Drinfeld's quantization conjecture, Etingof and Kazhdan constructed
a functor $Quant: \mathrm{LieBialg} \to \mathrm{Qgroup}$, which informally can be called
{\it universal quantization formula} or {\it quantum Baker--Campbell--Hausdorff formula}. They proved
that if one starts with a Lie bialgebra $L[[\hslash]]$ and applies the functor $Quant$ to it and then
$deQuant$, the resulting Lie bialgebra will be isomorphic to $L[[\hslash]]$. The same is true if one
starts with a quantum group $H$: $Quant(deQuant(H))$ will be isomorphic to $H.$

One of the uses of Lie groups-Lie algebras correspondence is the classification of semisimple Lie groups,
because the classification of semisimple Lie algebras is a much easier problem.
In the same way one can use Theorem \ref{KazEti} as an approach to the classification of quantum groups
over semisimple Lie algebras. This was done in the works \cite{SP}, \cite{SP2}.
The rest of the paper is dedicated to the exposition of the main results of that works.

\section{First steps of the classification}

Let $\mathfrak{g}$ be a simple Lie algebra over $\mathbb{C}$. We have seen that the classification of quantum groups over $\mathfrak{g}$ is equivalent to the classification of Lie bialgebra structures on $\mathfrak{g}[[\hslash]]:=\mathfrak{g} \otimes \mathbb{C}[[\hslash]]$. It is easy to see that any Lie bialgebra structure on $\mathfrak{g}[[\hslash]]$ gives rise to a Lie bialgebra structure on $\mathfrak{g}((\hslash)):=\mathfrak{g} \otimes \mathbb{C}((\hslash))$ and any Lie bialgebra structure on $\mathfrak{g}((\hslash))$  becomes a Lie biagebra structure on $\mathfrak{g}[[\hslash]]$ after a multiplication by an appropriate power of $\hslash$. Therefore, it is enough to classify Lie bialgebra structures on $\mathfrak{g}((\hslash))$.

Let us first look at the classification of Lie bialgebra structures on semisimple Lie algebras over an algebraically closed field $\mathbb{F}$ of characteristic zero. This classification was obtained by Belavin and Drinfeld \cite{BD}.  We will now give a brief outline of their results. Let $\delta$ be a Lie bialgebra structure on $\mathfrak{g}.$ First, one notices that the ``compatibility condition'' for $\delta$ is equivalent to the fact that $\delta$ is a cocycle. From the triviality of cohomology of simple Lie algebras we see that there exists $r \in \mathfrak{g} \otimes \mathfrak{g}$ such that $\delta=dr$.  The condition that $\delta$ is a Lie bialgebra structure can be rewritten in terms of $r$ : it turns out that after an appropriate scaling $r$ should satisfy the classical Yang-Baxter equation. There are two quite different cases, $r$ skew-symmetric  or non-skewsymmetric. In the first case there is no hope to obtain a meaningful classification. However, there is a lot of structure associated to a skew-symmetric $r$-matrix, this objects are intimately related to quasi-Frobenius Lie algebras \cite{BD}. In the second case Belavin and Drinfeld found the explicit formulas for $r$-matrices up to conjugation.
\begin{thm}
Let $\mathfrak{g}$ be a simple Lie algebra over an algebraically closed field of characteristic zero. Then any Lie bialgebra structure on $\mathfrak{g}$ is coboundary. Let $r$ be a corresponding $r$-matrix. If $r$ is not skew-symmetric then for some root decomposition we have
$$r = r_0 + \sum _ {\alpha > 0} e_{-\alpha}\otimes e_{\alpha} + \sum_{\alpha \in \mathrm{Span}( \Gamma_1)^{+}}\sum_{k \in \mathbb{N}}e_{-\alpha} \wedge e_{\tau^k(\alpha)} .$$
Here $\Gamma_1, \Gamma_2$ are the subsets of the set of simple roots, $\tau: \Gamma_1 \to \Gamma_2$ is isometric bijection, and for every $\alpha \in \Gamma_1$ there exists $k \in \mathbb{N}$ such that $\tau^k(\alpha) \in \Gamma_2 \setminus \Gamma_1$. The triple $(\Gamma_1, \Gamma_2,\tau)$ is called \emph{admissible}.
The tensor $r_0 \in \mathfrak{h} \otimes \mathfrak{h}$ must satisfy the following two conditions:

(1) $r_0+r_0^{21}=\sum t_k \otimes t_k$, where $t_k$ is an orthonormal basis of $\mathfrak{h}$,

(2) for any $\alpha \in \Gamma_1$ we have $(\tau(\alpha)\otimes \mathrm{id} + \mathrm{id} \otimes \alpha)r_0 = 0 $.
\end{thm}

It is worth noticing that there is an equivalent way to distinguish skew-symmetric and non-skewsymmetric $r$-matrices: in the first case the Drinfeld double $D(\mathfrak{g})$ is isomorphic to $\mathfrak{g} \otimes \mathbb{F}[\varepsilon]$, $\varepsilon^2=0$, in the second case $D(\mathfrak{g}) \simeq \mathfrak{g} \oplus \mathfrak{g}$, see \cite{St}.

 We want to obtain a version of Belavin-Drinfeld classification over the non-closed field $\mathbb{C}((\hslash))$. Let again $\mathfrak{g}$ be a simple Lie algebra over $\mathbb{C}$. First notice that we have a natural notion of equivalence for Lie bialgebras on $\mathfrak{g}((\hslash))$: $\delta_1 \sim \delta_2$ if and only if there exists $\lambda \in \mathbb{C}((\hslash))$ and $X \in G(\mathbb{C}((\hslash)))$ such that $\delta_1=\lambda \mathrm{Ad}_X \delta_2$. Here $G$ is an algebraic group associated to $\mathfrak{g}$

Any Lie bialgebra structure on $\mathfrak{g}((\hslash))$ can be lifted to $\mathfrak{g} \otimes \overline{\mathbb{C}((\hslash))}$. Over the algebraically closed field $\overline{\mathbb{C}((\hslash))}$ we have the Belavin-Drinfeld classification. Therefore, any Lie bialgebra structure on $\mathfrak{g}((\hslash))$ is given by an $r$-matrix of the form $\lambda \mathrm{Ad}_X r$, where $r$ is an $r$-matrix from the Belavin-Drinfeld list or a skew-symmetric $r$-matrix. One can prove that for a non-skew matrix up to equivalence $\lambda$ is either $1$ or $\sqrt{\hslash}$. Therefore, for any non-skew matrix from the Belavin-Drinfeld list there are two sets $H_{BD}^1(r_{BD})$ and $\overline{H}_{BD}^1 (r_{BD})$ of equivalence classes of $r$-matrices. $H_{BD}^1(r_{BD})$ parametrizes the equivalence classes of $r$-matrices of the form $\mathrm{Ad}_X r_{BD}$ that define a Lie bialgebra structure on $\mathfrak{g}((\hslash))$ and, respectively, $\overline{H}_{BD}^1 (r_{BD})$ parametrizes equivalence classes of matrices of the form $\sqrt{\hslash}\mathrm{Ad}_X r_{BD}$. We call $\overline{H}_{BD}^1(r_{BD})$ and $H_{BD}^1(r_{BD})$ the set of, respectively, twisted and non-twisted Belavin-Drinfeld cohomologies. Analogously for skew-symmetric $r$-matrix $r$ we define the Frobenius cohomology set $H_F^1(r)$.

There is an alternative way to see the difference between twisted and non-twisted Lie bialgebra structures. Let us look at the structure of the Drinfeld double $D(\mathfrak{g})$. It easily follows from methods developed in
\cite{MS} that there are three possible cases: $D(\mathfrak{g}((\hslash)))$
can be isomorphic to
$\mathfrak{g}((\hslash)) \oplus \mathfrak{g}((\hslash)),$
$\mathfrak{g}((\hslash))[\sqrt{\hslash}]$
or to
$\mathfrak{g}((\hslash))[\varepsilon]$, where $\varepsilon^2=0$.
These possibilities precisely correspond to the non-twisted,
twisted and skew cases respectively.

We have shown that all Lie bialgebra structures on $\mathfrak{g}$ fall into one of the three types: non-twisted, twisted or skew. In what follows we will examine each case in more detail.

\section{Non-twisted case}

We have defined $H_{BD}^1(r_{BD})$ as a set of equivalence classes of Lie bialgebra structures. However, there is an equivalent definition that appeals only to the inner structure of $\mathfrak{g}((\hslash))$. In what follows $G$ is an algebraic group that corresponds to $\mathfrak{g}$.
\begin{defn}
The \emph{centralizer} $C(r)$ of an $r$-matrix $r$ is the set of all $X \in G(\overline{\mathbb{C}((\hslash))})$ such that $\mathrm{Ad}_Xr=r$.
\end{defn}

\begin{defn}
$X \in G(\overline{\mathbb{C}((\hslash))})$ is called a \emph{non-twisted Belavin-Drinfeld cocycle} for $r_{BD}$ if for any $\sigma \in \Gal(\overline{\mathbb{C}((\hslash))}/\mathbb{C}((\hslash)))$ we have $X^{-1}\sigma(X) \in C(r_{BD})$. The set of non-twisted cocycles will be denoted by $Z(r_{BD})=Z(G, r_{BD})$.
\end{defn}

\begin{defn}
Two cocycles $X_1, X_2 \in Z(r_{BD})$ are called \emph{equivalent} if there exist $Q \in G(\mathbb{C}((\hslash)))$ and $C \in C(r_{BD})$ such that $X_1=QX_2C$.
\end{defn}

\begin{defn}
The set of equivalence classes of non-twisted cocycles is denoted by $H_{BD}^1(r_{BD})=H_{BD}^1(G, r_{BD})$ and is called the \emph{non-twisted Belavin-Drinfeld cohomology}.
\end{defn}

We were able to compute $H_{BD}^1$ for the algebras of $A-D$ series.
First let us make a small remark about $A_n$ case. In this case $\mathfrak{g}((\hslash))$ is naturally acted upon by the group $GL(n)$ and we can compute the cohomology with respect to conjugation by the elements of $GL(n)$ or $SL(n)$. To distinguish between these cases we write $H_{BD}^1(GL(n), r_{BD})$ and $H_{BD}^1(SL(n), r_{BD})$.

If $(\Gamma_1, \Gamma_2, \tau)$ is an admissible triple then the set
$\alpha, \tau(\alpha), ..., \tau^{k}(\alpha)$, where $\alpha \in \Gamma_1 \setminus \Gamma_2$ and $\tau^{k}(\alpha) \in \Gamma_2\setminus \Gamma_1$ will be called a string of $\tau$. The following table describes $H_{BD}^1$ for algebras of type $A-D$. The cohomology is called trivial if $|H_{BD}^1(r_{BD})|=1$.

\begin{table}[H]
\begin{tabular}{|c|c|c|c|}
\hline
Algebra                & Triple type                                                                                                                & \begin{tabular}[c]{@{}c@{}}$H_{BD}^1$ for \\an arbitrary field\end{tabular}  & $H_{BD}^1$ for $\mathbb{C}((\hslash))$ \\ \hline
$A_n$                  &                                                                                                                            & trivial ($GL(n)$ case)                        &                                      \\ \hline
$B_n$                  &                                                                                                                            & trivial                         &                                      \\ \hline
$C_n$                  &                                                                                                                            & trivial                         &                                      \\ \hline
\multirow{2}{*}{$D_n$} & \begin{tabular}[c]{@{}c@{}}there exists a \\ string of $\tau$ that\\ contains $\alpha_{n-1}$\\ and $\alpha_n$\end{tabular} & $F^*/(F^*)^2$ & 2 elements                           \\ \cline{2-4}
                       & \begin{tabular}[c]{@{}c@{}}$\alpha_{n-1}$ and $\alpha_n$\\ do not belong to \\ the same string of $\tau$\end{tabular}      & trivial                         & \multicolumn{1}{l|}{}                \\ \hline
\end{tabular}
\end{table}
\begin{rem} In this paper $\alpha_n,\alpha_{n-1}$ are the branchendpoints in the Dynkin diagram for $D_n$
\end{rem}
\begin{rem}
One can similarly define Belavin-Drinfeld cohomologies over an arbitrary field $F$ as a tool  to understand Lie bialgebra structures on $\mathfrak{g}(F)$.

\end{rem}

The result for $H_{BD}^1(SL(n), r_{BD})$ is more interesting. Let $\alpha_{i_1}, \ldots , \alpha_{i_k}$ be a string of $\tau$, $\tau (\alpha_{i_p})=\alpha_{i_{p+1}}$.
If $\tau (\alpha_{i_p})$ is not defined, then anyway we define the corresponding string, which consists of one element $\{\alpha_{i_p}\}$ only.

For any string $S=\{\alpha_{i_1}, \ldots , \alpha_{i_k}\}$ of $\tau$, we define the weight of $S$ by $w_S=\sum_p i_p$. Moreover, for any Belavin-Drinfeld triple we will also formally consider the string $\{\alpha_n\}$ with weight $n$. Let $t_1, \ldots , t_n$ be the ends of the strings with weights $w_1, \ldots , w_n$. We note that some indices in  $w_1, \ldots , w_n$ are missing unless $\Gamma_1$ is an empty set and $w_n=n$ is always present. Let $N=GCD(w_1, \ldots , w_n)$.
\begin{thm}
The number of elements of $H^1_{BD}(SL(n), r)$ is $N$. Each cohomology class contains a diagonal matrix $D=A_1A_2$,
where $A_2\in C(GL(n), r)$ and
$A_1\in\mathrm{diag}(n, \mathbb{C}((\hslash)))$.
Two such diagonal matrices $D_1=A_1A_2$  and $D_2=B_1B_2$ are contained in the same class of $H_{BD}^1(SL(n), r)$
if and only if $\det(A_1)=\det(B_1)$
in ${\mathbb{C}((\hslash))}^*/(\mathbb{C}((\hslash))^*)^N $.
\end{thm}
\section{Twisted case}
As in non-twisted case there is a way to define $\overline{H}_{BD}^1$ without mentioning Lie bialgebra structures.
\begin{thm}
$a \mathrm{Ad}_X r_{BD}$ defines a Lie bialgebra structure on $\mathfrak{g}(\mathbb{C}((\hslash)))$ if and only if $X$ is a non-twisted cocycle for the field $\mathbb{C}((\hslash))[\sqrt{\hslash}]$ and $\mathrm{Ad}_{X^{-1}\sigma_0(X)}r_{BD}=r_{BD}^{21}$. Here $\sigma_0$ is the non-trivial element of the group $\Gal(\mathbb{C}((\hslash))[\sqrt{\hslash}]/\mathbb{C}((\hslash)))$.
 \end{thm}
To deal with the condition  $\mathrm{Ad}_{X^{-1}\sigma_0(X)}r_{BD}=r_{BD}^{21}$ we classified all triples $(\Gamma_1, \Gamma_2, \tau)$ such that $r_{BD}^{21}$ and $r_{BD}$ are conjugate. In each case we found a suitable $S \in G(F)$ such that $r_{BD}^{21}=\mathrm{Ad}_S r_{BD}$. Then we can define Belavin-Drinfeld cocycles and cohomologies similar to the non-twisted case. In all cases $S^2=\pm 1$.
\begin{defn}
 $X \in G(\overline{\mathbb{C}((\hslash))})$ is called a \emph{Belavin-Drinfeld twisted cocycle} if for any $\sigma \in \Gal(\overline{\mathbb{C}((\hslash))}/\mathbb{C}((\hslash))[\sqrt{\hslash}])$ we have $X^{-1}\sigma(X) \in C(r_{BD})$ and $SX^{-1}\sigma_0(X) \in C(r_{BD})$. The set of Belavin-Drinfeld twisted cocycles is denoted by $\overline{Z}(r_{BD})=\overline{Z}(G, r_{BD})$.
\end{defn}
\begin{defn}
 Two twisted cocycles $X_1, X_2$ are called \emph{equivalent} if there exist $Q \in G(\mathbb{C}((\hslash)))$ and $C \in C(r_{BD})$ such that $X_1=QX_2C$. The set of equivalence classes of twisted cocycles is called the \emph{twisted Belavin-Drinfeld cohomology} and is denoted by $\overline{H}_{BD}^1(r_{BD})=\overline{H}_{BD}^1(G, r_{BD})$.
 \end{defn}
\begin{table}[h]
\begin{tabular}{|c|cc|c|}
\hline
Algebra                & \multicolumn{2}{c|}{Triple type}                                                                                                                                                                                                                                                                                                                                                                                                             & $\overline{H}_{BD}^1$ for $\mathbb{C}((\hslash))$ \\ \hline
\multirow{2}{*}{$A_n$} & \multicolumn{2}{c|}{\begin{tabular}[c]{@{}c@{}}$s\tau=\tau^{-1}s$, where $s$ is the non-trivial\\ automorphism of the Dynkin diagram\end{tabular}}                                                                                                                                                                                                                                                                                                       & one element                                     \\ \cline{2-4}
                       & \multicolumn{2}{c|}{other}                                                                                                                                                                                                                                                                                                                                                                                                                   & empty                                           \\ \hline
\multirow{2}{*}{$B_n$} & \multicolumn{2}{c}{Drinfeld-Jimbo}                                                                                                                                                                                                                                                                                                                                                                                                           & one element                                     \\ \cline{2-4}
                       & \multicolumn{2}{c|}{not DJ}                                                                                                                                                                                                                                                                                                                                                                                                                  & empty                                           \\ \hline
\multirow{2}{*}{$C_n$} & \multicolumn{2}{c|}{Drinfeld-Jimbo}                                                                                                                                                                                                                                                                                                                                                                                                          & one element                                     \\ \cline{2-4}
                       & \multicolumn{2}{c|}{not DJ}                                                                                                                                                                                                                                                                                                                                                                                                                  & empty                                           \\ \hline
\multirow{5}{*}{$D_n$} & \multicolumn{1}{c|}{\multirow{2}{*}{even $n$}} & Drinfeld-Jimbo                                                                                                                                                                                                                                                                                                                                                                              & one element                                     \\ \cline{3-4}
                       & \multicolumn{1}{c|}{}                          & not DJ                                                                                                                                                                                                                                                                                                                                                                                      & empty                                           \\ \cline{2-4}
                       & \multicolumn{1}{c|}{\multirow{3}{*}{odd $n$}}  & \begin{tabular}[c]{@{}c@{}}$\Gamma_1=\{\alpha_{n-1}\}$\\ $\tau(\alpha_{n-1})=\alpha_n$ ;\\ $\Gamma_1=\{\alpha_n\}$\\ $\tau(\alpha_n)=\alpha_{n-1}$ ;\\ $\Gamma_1=(\alpha_{n-1}, \alpha_k)$, $k \neq n$\\ $\tau(\alpha_{n-1})=\alpha_k, \tau(\alpha(k))=\alpha_n$;\\ $\Gamma_1=(\alpha_{n}, \alpha_k)$, $k \neq n-1$\\ $\tau(\alpha_{n}=\alpha_k), \tau(\alpha_k)=\alpha_{n-1}$\end{tabular} & two elements                                    \\ \cline{3-4}
                       & \multicolumn{1}{c|}{}                          & Drinfeld-Jimbo                                                                                                                                                                                                                                                                                                                                                                              & one element                                     \\ \cline{3-4}
                       & \multicolumn{1}{c|}{}                          & not DJ                                                                                                                                                                                                                                                                                                                                                                                      & empty                                           \\ \hline
\end{tabular}
\end{table}
Here cohomologies for $\mathfrak{sl}_n$ are considered with respect to the group $GL(n)$. For the results for $A_n$ over arbitrary field see \cite{PS}.

\section{Skew-symmetric case}
Following the pattern in  \cite{BD}, it can be easily proved that the classification of Lie bialgebra structures related to skew (triangular) $r$-matrices on $\mathfrak{g}((\hslash))$ is equivalent to the classification of quasi-Frobenius Lie subalgebras of $\mathfrak{g}((\hslash))$. This can be used to prove that if $r$ is skew-symmetric  then $r$ has to be defined over $\mathbb{C}((\hslash))$. However, different $r$-matrices defined over $\mathbb{C}((\hslash))$ can be conjugate over $\overline{\mathbb{C}((\hslash))}$. We can define Frobenius cohomology similarly  to Belavin-Drinfeld cohomology. We call two $r$-matrices equivalent if there exists $a \in {\mathbb{C}((\hslash))}$, $X \in G(\mathbb{C}((\hslash)))$ such that $r_1= a \mathrm{Ad}_X  r_2$. If $r$ defines a Lie bialgebra structure on $\mathfrak{g}((\hslash))$ then we define the Frobenius cohomology set $H_F^1(r)$ to be the set of equivalence classes of $r$-matrices that are conjugate to $r$ over $\overline{\mathbb{C}((\hslash))}$. Even though we do not have a classification of skew $r$-matrices even over an algebraically closed field, this cohomology can be computed in a way similar to Belavin-Drinfeld case.
\begin{defn}
The \emph{centralizer} $C(r)$ of an $r$-matrix $r$ is the set of all $X \in G(\overline{\mathbb{C}((\hslash))})$ such that $\mathrm{Ad}_Xr=r$.
\end{defn}

\begin{defn}
$X \in G(\overline{\mathbb{C}((\hslash))})$ is called a \emph{non-twisted Frobenius cocycle} for $r$ if for any $\sigma \in \Gal(\overline{\mathbb{C}((\hslash))}/\mathbb{C}((\hslash)))$ we have $X^{-1}\sigma(X) \in C(r)$. The set of non-twisted cocycles will be denoted by $Z_{F}(r)=Z_{F}(G, r)$.
\end{defn}

\begin{defn}
Two cocycles $X_1, X_2 \in Z_F(r)$ are called \emph{equivalent} if there exists $Q \in G(\mathbb{C}((\hslash)))$ and $C \in C(r)$ such that $X_1=QX_2C$.
\end{defn}

\begin{defn}
The set of equivalence classes of Frobenius cocycles is denoted by $H_F^1(r)=H_F^1(G, r)$ and is called the \emph{Frobenius cohomology}.
\end{defn}
\begin{ex} Let $r_J$ be the  Jordan $r$-matrix, i.e. $r_J=E \wedge H$. Then $H_F^1(r_J)$ is trivial. Here $\{E,F,H\}$ is the standard basis in $\frak{s}\frak{l}_2$
\end{ex}

\section{Historical Remarks}
Quantum groups (as in Definition \ref{Qgr}) were defined by Drinfeld in his talk at the International Congress of Mathematicians  in Berkeley, 1986. Relations between quantum groups and quantum algebras (quantization and dequantization functors, quantum Baker-Campbell-Hausdorff formula) were obtained by Etingof and Kazhdan in a series of papers \cite{EK1},\cite{EK2}.

The first example of a quantum group of the non-twisted type is due to Kulish and Reshetikhin \cite{KR}.
Generalizations for all simple Lie algebras were obtained by Drinfeld and Jimbo \cite{Dr},\cite{Ji}, where they
found quantum groups which quantize Lie bialgebra structures on $\mathfrak{g}$ defined by $\Gamma_1=\Gamma_2=\emptyset$.

Further classes of Lie bialgebra structures on $\frak{g}$, related to certain triples $(\Gamma_1,\Gamma_2,\tau)$, were quantized by Kulish and Mudrov in \cite{KM}.

Finally, Etingof, Schiffman, and Schedler quantized all Lie bialgebra structures
defined by all admissible triples  $(\Gamma_1,\Gamma_2, \tau)$ \cite{ESS}.

There are no explicit formulas for quantum groups related to the twisted Belavin-Drinfeld cohomologies.

Construction of quantum groups of the skew-symmetric  type appeared in the work of Drinfeld \cite{DR1} by means of a certain twisting element $F$.
The first explicit formula for $F$ is due to Coll, Gerstenhaber, and Giaquinto \cite{CGG}. This formula was used by Kulish and Stolin to quantize a certain nonstandard Lie bialgebra structure on the polynomial Lie algebra $\frak{sl}_2[u]$.

This paper is dedicated to Petr P. Kulish on the occasion of his 70-years jubilee.  The authors are thankful for valuable remarks to G.~Rozenblum who joins the congratulations.

\end{sloppy}
\end{document}